\documentclass{ifacconf}
\usepackage{graphicx}      
\usepackage{natbib}        

\usepackage{epsfig} 
\usepackage{amsmath} 
\usepackage{amssymb}  
\usepackage{float}
\begin{document}
\begin{frontmatter}

\title{A Smart Sensor for the Heat Flow of an Extruder
}

\author[First]{K. Schwarzinger} 
\author[First]{K. Schlacher}

\address[First]{Johannes Kepler University Linz, Institute of Automatic Control and Control Systems Technology, Altenberger Straße 69, A-4040 Linz, Austria, (e-mail: $ \{ $kevin.schwarzinger, kurt.schlacher$ \} $@jku.at)}

\begin{abstract}                
	Two mathematical models are derived to describe the thermal dynamics of a plastics production machine. The Smart Sensor design is based on the finite volume model. The accuracy of the derived finite volume model is validated by a finite element model. Those models are validated by measurements. A Smart Sensor for the heat flow between the extruder cylinder and screw conveyor is designed. Such a configuration has the advantage to decouple the thermal control problem from the specific granulate.
\end{abstract}


\end{frontmatter}

\section{Introduction}
Extruders represent essential parts in the plastic production process as they convert the granulate to melt. Its control performance has a great impact on the work machines and therefore on the quality of the products. There exists several contributions which are devoted to mathematical modelling, simulation and optimization of the production process, see, e.$ \, $g., \cite{KOCHHAR1977177}, \cite{tadmor2013principles}, \cite{acur1982numerical}, \cite{Zheng2012} and \cite{Ravi2010ModelExtruder}.\\
Here, a precise mathematical model is constructed, so that the system order is as small as possible. The model is derived following the method of finite volumes (FV). Its simulation outcome is confronted to the numeric results from a precise finite element (FE) model, generated by the simulation package \textit{FEniCS} (described in \cite{LoggOlgaardEtAl2012a} and \cite{AlnaesBlechta2015a}). A similar FE approach is described by \cite{belavy2018fem}. These models are the basis for the design of a Smart Sensor which estimates the heat flow between the screw conveyor and the cylinder. The cylinder is actuated by electrical heating tapes and equipped with thermocouples.\\ 
The paper is organized as follows: The extruder setup is discussed in Section \ref{ch:extruderM}. In Section \ref{ch:model} the FV as well as the FE approach are described and the simulation results of both models are compared with measurements from an industrial device.
The design of the Smart Sensor is shown in Section \ref{ch:HeatFlowO} together with experimental results by different extrusion experiments.

\section{The extruder setup}\label{ch:extruderM}
\begin{figure}[tbph]
	\centering
	\includegraphics[scale=1.0]{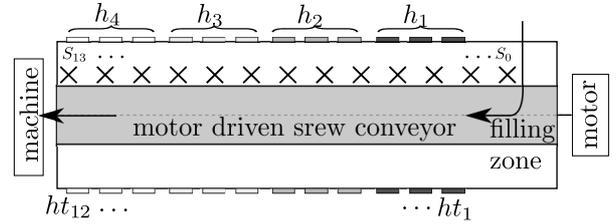}
	\caption{A typical extruder configuration.}
	\label{fig:extruder}
\end{figure}

The setup of a typical extruder is shown in  Fig. \ref{fig:extruder}, it consists of a cylinder with length $ L $, a inner diameter $ d_1 $ and an outer diameter $ d_2 $. The screw conveyor is driven by a motor. The granulate is admixed by a dosing device in the filling zone. $ N_{th} $ heating tapes are mounted around the extruder cylinder. A common extruder setup comprises $ N_{h} $ heating zones, which are formed by several heating tapes. The heating tapes are connected in parallel. The number of heating tapes per heating zone is given by $ N_{th,h} = N_{th}/N_{h} $. The state of the extruder is given by the temperature distribution $ T(t, \overline{\mathbf{x}}) $ at time $ t $ and the location coordinate $ \overline{\mathbf{x}} $. The thermocouples ($S_{0} $ - $ S_{13} $ in Fig. \ref{fig:extruder}) provide us with the temperature at the place of installation.
\\
It is a common approach to combine a heating zone and a centred thermocouple from the heating zone (one of $S_{0} $ - $ S_{13} $) inside the cylinder with a PID-law. This simple approach does not require a mathematical model. Since the coupling between the heating zones is neglected the performance of the controlled system is limited. Here we are interested in a mathematical model of the whole cylinder for advanced control strategies where the coupling of the heating zones is precisely taken into account. The classic PID approach does neither consider an interaction between the heating zones nor production-related desires concerning local temperature restrictions, which a controller based on the computed model is indeed capable of.\\
The heating tapes are actuated by a switching device. The input for this device is a number of pulses in a fixed time interval. Since there exists no simple relation between the number of pulses and the heat flow injected into the cylinder, a subordinate controller is used for the control of the temperature given by a thermocouple with output $ T_{h} $ mounted between cylinder and the centred heating tape. These control loops are fast and allow us to neglect the coupling of the heating zones. This should not be mixed up with the approach described above where the sensors are mounted far from the heating zones.
One should keep in mind, that this configuration does not allow active cooling, but the measured $ T_{h} $ can be considered as input to the mathematical model of the cylinder. Even if the system equations are linear, the control problem is highly non-linear because of the restriction (no active cooling) in the control input.

\section{Mathematical modelling / Validation}\label{ch:model}
\subsection{Finite Volume Approach}\label{ch:FV}
The extruder to be modelled consists of rotationally symmetrical heating tapes installed around the cylinder. Therefore, the model can be described as a rotationally symmetrical problem. Due to the rotational symmetry assumption, the position is determined by two coordinates $ \overline{\mathbf{x}} = \begin{pmatrix}x& r\end{pmatrix}^T  $. The governing equations for our problem are given by the Fourier heat equation in cylindrical coordinates 
\begin{align}\label{eq:heateq1}
	\begin{aligned}
		\rho(\overline{\mathbf{x}}) c_{p}(\overline{\mathbf{x}},T) \frac{\partial T(t, \overline{\mathbf{x}})}{\partial t}=
		\frac{1}{r} \frac{\partial}{\partial r}\left(\lambda(\overline{\mathbf{x}},T)r \frac{\partial T(t, \overline{\mathbf{x}})}{\partial r}\right)\\
		+\frac{\partial}{\partial x}\left(\lambda(\overline{\mathbf{x}},T)\frac{\partial T(t, \overline{\mathbf{x}})}{\partial x}\right).
	\end{aligned}
\end{align}
The parameters $ \rho$, $c_{p} $ and $ \lambda $ (see Tab. \ref{tb:parameter}) are assumed to be constant inside each individual FV element.
\begin{table}[H]
	\caption{Notations for the model.}
	\label{tb:parameter}
	\begin{center}
		\begin{tabular}{ccc}
			\hline
			Name & Symbol & Unit \\
			\hline
			time & $ t $ & $ \mathrm{s} $\\
			temperature & $ T $ & $ \mathrm{K} $\\
			density & $ \rho $ &  $ \frac{\mathrm{kg}}{\mathrm{m}^3} $\\
			specific heat capacity & $ c_p $ &  $ \frac{\mathrm{J}}{\mathrm{kg}\mathrm{K}} $\\
			thermal conductivity & $ \lambda $ &  $ \frac{\mathrm{W}}{\mathrm{m}\mathrm{K}} $\\
			\hline
		\end{tabular}
	\end{center}
\end{table}
Due to the rotationally symmetry assumptions, a FV element is described by a ring element with length $ \Delta x_i $ and the radial thickness $ \Delta r_i $.
\begin{figure}[H]
	\centering
	\includegraphics[scale=1.0]{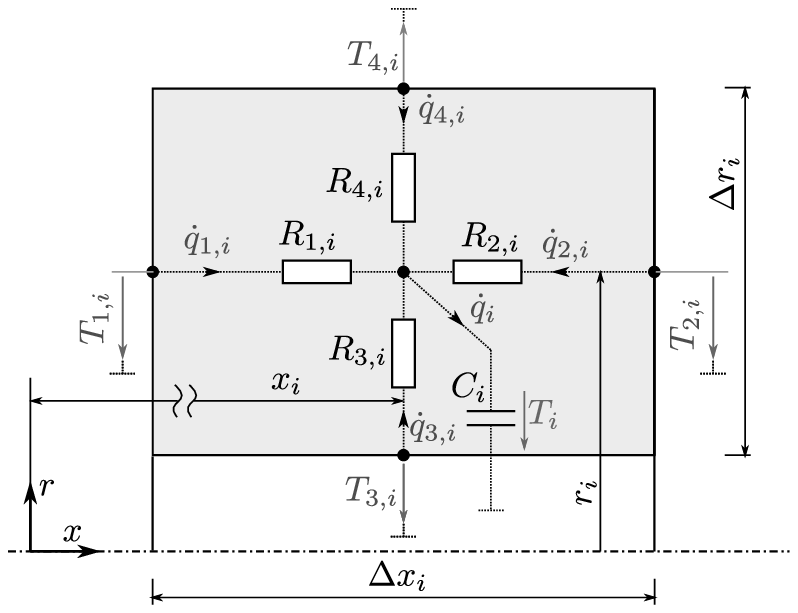}
	\caption{Equivalent circuit for a typical extruder configuration.}
	\label{fig:eleESB}
\end{figure}
The FV elements are connected by ports, so that a number of elements $ n_a $ in the axial and  $ n_r $ in the radial direction arises, resulting in a total number of $ n = n_a n_r $ elements. The number of elements correspond to the resulting system order. The surfaces of the $ i $-th FV are named $ A_i $ and the volume $ V_i $.
According to Fig. \ref{fig:eleESB} the following abbreviations $T(t,r, x-\Delta x) = T_{1,i}$, $T(t,r, x+\Delta x) = T_{2,i}$, $T(t,x, r-\Delta r) = T_{3,i}$ and $T(t,x, r+\Delta r) = T_{4,i}$ are used. The size of the elements are denoted by $ \Delta r_i$ and $\Delta x_i$, this allows us to combine FV elements of different sizes. Different volumes sizes enables a simple adaption of the elements to the cylinder including heating tapes and thermocouples. The following equations
\begin{align}\label{eq:diskheatequ_RC_FV}
	\begin{aligned}
		C_i \frac{\partial T_i}{\partial t}=&\frac{1}{R_{1,i}}\left(T_{1,i}-T_{i}\right)+\frac{1}{R_{2,i}}\left(T_{2,i}-T_{i}\right)\\
		+&\frac{1}{R_{3,i}}\left(T_{3,i}-T_{i}\right)+\frac{1}{R_{4,i}}\left(T_{4,i}-T_{i}\right),
	\end{aligned}
\end{align}
with the quantities 
\begin{align}\label{eq:RC_val_1}
	\begin{aligned}
		&C_i=2 r_i \Delta r_i \pi \Delta x_i \rho_i c_{p,i}=V_i\rho_i c_{p,i}\\
		&R_{[1,2],i}=\frac{\Delta x_i}{2r_i \Delta r_i \pi\lambda_i}=\frac{\Delta x_i}{ A_{[1,2],i}\lambda_i}\\
		&R_{3,i}=\frac{\Delta r_i}{2\pi(r_i+\frac{\Delta r_i}{2})\Delta x_i\lambda_i}=\frac{\Delta r_i}{A_{3,i}\lambda_i}\\
		&R_{4,i}=\frac{\Delta r_i}{2\pi(r_i-\frac{\Delta r_i}{2})\Delta x_i\lambda_i}=\frac{\Delta r_i}{A_{4,i}\lambda_i},
	\end{aligned}
\end{align}
represent the spatially discretised model according to Fig. \ref{fig:eleESB}. 
It is straightforward by combining \eqref{eq:diskheatequ_RC_FV} and \eqref{eq:RC_val_1} to derive
\begin{align}\label{eq:diskheatequ}
	\begin{aligned}
		\rho_i c_{p,i}& \frac{\partial T_i}{\partial t}=\lambda_i \frac{T_{1,i}-2T_i+T_{2,i}}{\Delta x_i^2}\\
		+&\lambda_i \frac{T_{3,i}-2T_i+T_{4,i}}{\Delta r_i^2} + \frac{\lambda_i}{r_i}\frac{T_{3,i}-T_{4,i}}{2 \Delta r_i}.
	\end{aligned}
\end{align} 
Replacing the central difference quotients by the corresponding first and second order derivatives yields
\begin{align}\label{eq:limit}
	\begin{aligned}
		& \rho_i c_{p,i} \frac{\partial T_i}{\partial t}=\lambda_i\left(\frac{\partial^2 T_i}{\partial x^2} + \frac{1}{r_i}\frac{\partial T_i}{\partial r}+\frac{\partial^2 T_i}{\partial r^2}\right),
	\end{aligned}
\end{align}
which coincided with \eqref{eq:heateq1}, confirming that the heat equation approximately corresponds to the model of an electrical network which consists of a capacitor and resistors.\\\\
The input of the network are temperature sources or heat flow sources. Each FV element has four ports. The port variables are $ T_{s,i} $ and $ \dot{q}_{s,i} $, which can be combined with the neighbour or with the environment. In addition, the port variables are connected by one linear equation, see, Fig. 2, for the port variables and the temperature $ T_i $ of the FV element. 
In particular the interconnection with a neighbour element $\hat{i}$, where index $ s $ and $ \hat{s} $ mark nodes of the $ i $-th and $\hat{i}$-th FV element, respectively, is given by
\begin{align}\label{eq:connectNode}
	\begin{aligned}
		T_{s,i}=T_{\hat{s},\hat{i}}&&\dot{q}_{s,i}=-\dot{q}_{\hat{s},\hat{i}}
	\end{aligned}\,.
\end{align}
This method has to be applied iteratively to generate a set of additional equations. To complete the model construction, boundary conditions (BC) have to be defined. Simple connections to the environment are given by
\begin{subequations}
	\begin{align}
		T_{x,i}(t,\overline{\mathbf{x}}) &= T_{N}(t,\overline{\mathbf{x}}) \label{eq:DRB}\\
		\dot{q}_{x,i}(t,\overline{\mathbf{x}}) &= \dot{q}_{N}(t,\overline{\mathbf{x}}), \label{eq:NRB}
	\end{align}
\end{subequations}
where the quantities $ \dot{q}_{x,i} $ and $ T_{x,i} $ follow from the network equations. These are the well known Dirichlet \eqref{eq:DRB} and Neumann \eqref{eq:NRB} boundary conditions. Other boundary conditions are
\begin{subequations}
	\begin{align}
		\dot{q}_{x,i}(t,\overline{\mathbf{x}}) &= \alpha_i(\overline{\mathbf{x}})\left(T_{x,i}(t,\overline{\mathbf{x}})-T_{0}(t,\overline{\mathbf{x}})\right)  \label{eq:RRB} \\
		\dot{q}_{x,i}(t,\overline{\mathbf{x}}) &= \alpha_i(\overline{\mathbf{x}})\left(T_{x,i}(t,\overline{\mathbf{x}})-T_{0}(t,\overline{\mathbf{x}})\right) \nonumber\\
		& +\alpha_{ht,i}(\overline{\mathbf{x}})\left(T_{x,i}(t,\overline{\mathbf{x}})-T_{h}(t,\overline{\mathbf{x}})\right)  \label{eq:HTRB} \\
		\dot{q}_{x,i}(t,\overline{\mathbf{x}}) &= \alpha T_{x,i}^4(t,\overline{\mathbf{x}}) .\label{eq:STRB}
	\end{align}
\end{subequations}
The connection with the environment is established by \eqref{eq:RRB}, whereas $ T_{0} $ is the ambient temperature and $ \alpha_{i} $ is the heat transfer coefficient between cylinder and environment. The heating tapes are defined by \eqref{eq:HTRB}, whereas $ \alpha_{ht,i} $ is the heat transfer coefficient between cylinder and the heating tape and can be interpreted as depicted in Fig. \ref{fig:zones} on the left. The presented discretization allows us to take all types of boundary conditions into account.\\
If we combine all the equations of our network \eqref{eq:diskheatequ_RC_FV} with \eqref{eq:connectNode}, \eqref{eq:DRB} - \eqref{eq:NRB} and \eqref{eq:RRB} - \eqref{eq:STRB}, we derive a linear system of the form
\begin{align}\label{eq:eqset}
	\begin{aligned}
		\dot{\mathbf{x}} &= \mathbf{F}\mathbf{x}+\mathbf{H}\mathbf{d} + \mathbf{G}\mathbf{u}\\
		\mathbf{0}&= \tilde{\mathbf{F}}\mathbf{x}+\tilde{\mathbf{H}}\mathbf{d} + \tilde{\mathbf{G}}\mathbf{u}\\
	\end{aligned}
\end{align}
with an invertible block diagonal matrix $ \tilde{\mathbf{H}} $ and the eliminable quantities from the network equations $ \mathbf{d} $. Therefore, it is straight forward to derive the linear state space model
\begin{align}\label{eq:result}
	\begin{aligned}
		\dot{\mathbf{x}} &= \mathbf{A}\mathbf{x}+\mathbf{B}\mathbf{u}=\begin{bmatrix}
			\mathbf{A}_{11} & \mathbf{A}_{12} & \mathbf{0}\\
			\mathbf{A}_{21} & \mathbf{A}_{22} & \mathbf{A}_{23}\\
			\mathbf{0} & \mathbf{A}_{32} & \mathbf{A}_{33}
		\end{bmatrix}\begin{bmatrix}
		\mathbf{x}_c\\
		\mathbf{x}_s\\
		\mathbf{x}_{sc}
	\end{bmatrix}+\begin{bmatrix}
	\mathbf{B}_c\\
	\mathbf{B}_s\\
	\mathbf{B}_{sc}
\end{bmatrix}\mathbf{u}\\
		\mathbf{y} &= \mathbf{C}\mathbf{x}= \begin{bmatrix}
			\mathbf{C}_1 & \mathbf{0} & \mathbf{0}
		\end{bmatrix}\mathbf{x}\,.
	\end{aligned}
\end{align}
The state is given by the temperatures $ T_{i} $ of the FV elements and input by the temperature sources
\begin{align}\label{eq:resultux}
	\begin{aligned}
		\mathbf{u} &= \begin{bmatrix} T_{h,1} & \cdots & T_{h,N_{h}} & T_0 \end{bmatrix}^T\\
		\mathbf{x} &= \begin{bmatrix} \mathbf{x}_c^T &  \mathbf{x}_s^T & \mathbf{x}_{sc}^T \end{bmatrix}^T = \begin{bmatrix} T_1 & \cdots & T_{n\cdot m} \end{bmatrix}^T\,.
	\end{aligned}
\end{align}
The state vector $ \mathbf{x} $ consists of all FV temperatures $ T_{i}$ and can be divided into the temperatures belonging to the cylinder zone $ \mathbf{x}_c $, to the screw conveyor zone $ \mathbf{x}_c $ and the screw core zone $ \mathbf{x}_{sc} $. A description of the introduced zones follows below. The input vector $ \mathbf{u} $ contains the ambient temperature $ T_{0}$ and all temperatures sources\footnote{In the considered setup there are no heat flow sources.}  $ T_{h,j}$. The output $ \mathbf{y} $ includes the temperatures of the points where the sensors are placed (the points $ S_i $ in Fig. \ref{fig:extruder}).\\
\begin{figure}[H]
	\centering
	\includegraphics[scale=1.0]{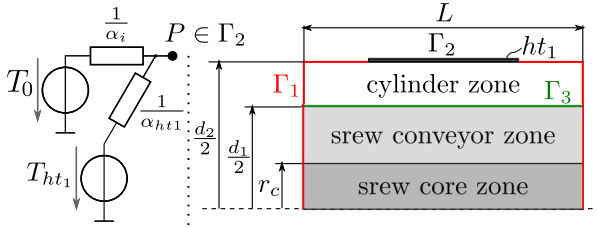}
	\caption{BC realisation for nodes $ P\in\Gamma_2 $ on the left, boundary decomposition in sets $ \Gamma_1 $ to $ \Gamma_2 $ and the zone definition at the right.}
	\label{fig:zones}
\end{figure}
The right side of Fig. \ref{fig:zones} illustrates the considered modelling area, which is described by an outer diameter $ d_2 $ and the extruder length $ L $. Three radially limited material related areas are defined: from the centre axis to the core of the screw conveyor $ r_c $ - \textit{screw core zone} ($ c_{p,sc} $, $ \lambda_{sc} $ and $ \rho_{sc} $), the area form the core of the screw conveyor to the inner radius of the extruder cylinder $ d_1/2 $ -  \textit{screw conveyor zone} ($ c_{p,s} $, $ \lambda_{s} $ and $ \rho_{s} $) and from $ d_1/2 $ to the outer radius of the extruder $ d_2/2 $ - \textit{cylinder zone} ($ c_{p,c} $, $ \lambda_{c} $ and $ \rho_{c} $).
In Fig. \ref{fig:zones} two sets $ \Gamma_1 $ to $ \Gamma_2 $, consisting of all points that are located on the border area $ \Gamma = \Gamma_1 \cup \Gamma_2 $, can be seen.  The set $ \Gamma_2 $ includes open nodes that are directly located under the heating tapes and $ \Gamma_1 $ includes nodes that are in direct contact with the environment. The most inner FV elements are disks not rings.\\
To monitor the heat flow across the border of the extruder and the heat flow injected by the heating tapes, we can extend the output with additional equations. A big advantage this model offers is that the operator gets information about the temperature distribution as well as the heat flows inside or out of the extruder. The mutual influence of the heating tapes is also taken into account.

\begin{figure*}[tbph]
	\centering
	\includegraphics[scale=0.7]{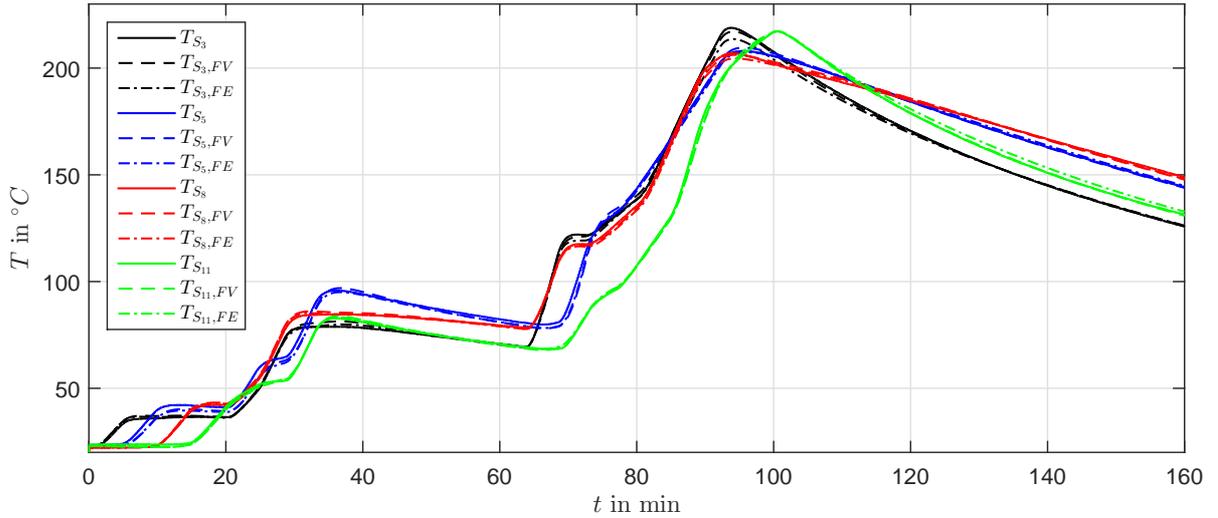}
	\caption{Comparison between the measurement and the FV- and FE - model simulation results.}
	\label{fig:Presult}
\end{figure*}
\subsection{Finite Element Approach}\label{ch:FE}
A model created using the finite element (FE) approach is used to derive an accurate reference solution to find an FV model \eqref{eq:result} with the lowest possible system order that is comparably accurate. The FE model comes with a much higher system order resulting in a higher resolution and accuracy of the computed temperature distribution. In order to derive a FE model, a variational formulation of the partial differential equations is needed. The locational functions describing $ c_p(\overline{\mathbf{x}})$, $ \lambda(\overline{\mathbf{x}}) $ and $ \rho(\overline{\mathbf{x}}) $ have to be defined so that the different zones, illustrated in Fig. \ref{fig:zones}, are taken into account.
A FE program like \textit{FEniCS} ,see \cite{LoggOlgaardEtAl2012a}, \cite{AlnaesBlechta2015a} can only solve static problems.
By help of the Backward Euler formula
\begin{align}\label{eq:beuler}
	\frac{\partial T(t, \overline{\mathbf{x}})}{\partial t} = \frac{T(k)-T(k-1)}{\Delta t}
\end{align}
the dynamic problem is replaced by a sequence of static ones with the time step $ \Delta t $. The variational formulation with the test functions $ v\in V $ is given by 
\begin{align} \label{eq:var6}
	\begin{aligned}
		&\int_{\Omega}\rho c_{p}T(k-1)v\mathrm{~d}A=-\Delta t\lambda\int_{\Omega}\frac{1}{r}\frac{\partial T(k)}{\partial r}v\mathrm{~d}A~+
		\\&\int_{\Omega}\left[\rho c_{p}T(k)v+\Delta t\lambda\left(\frac{\partial T(k)}{\partial r}\frac{\partial v}{\partial r}+\frac{\partial T(k)}{\partial x}\frac{\partial v}{\partial x}\right)\right]\mathrm{d}A
		\\&-\Delta t\lambda\int_{\Gamma_{1}}\alpha_{\gamma_{1}}(T(k)-T_{u})v\mathrm{~d}s
		\\&-\Delta t\lambda\int_{\Gamma_{2}}\left(\alpha_{\gamma_{2}}(T(k)-T_{u})+\alpha_{ht,i}(T(k)-T_{h,i})\right)v\mathrm{~d}s,
	\end{aligned}
\end{align} 
where $ \alpha_{\gamma_1} $ and $ \alpha_{\gamma_2} $ denote the surface independent heat transfer coefficients between the extruder and its surrounding and $ \alpha_{ht,i} $ describe the surface independent heat transfer coefficients between the $ i $-th heating tape and the extruder. The space of allowed test functions $ V=\left\{v(\overline{\mathbf{x}}) \in H^{1}(\Omega): v(\overline{\mathbf{x}})=0, \forall \overline{\mathbf{x}} \in \Gamma_\alpha \right\}  $ is the Sobolev space $  H^{1} $ where  $ \Omega $ is the area of integration and $ \Gamma_\alpha $ its boundaries on which a BC in the form of \eqref{eq:DRB} is going to be implemented. Because of space limitations the derivation of \eqref{eq:var6} is not discussed in detail.
\subsection{Validation Through Heat Up Experiments} \label{ch:heatupT}
\begin{figure}[htbp]
	\centering
	\includegraphics[scale=1]{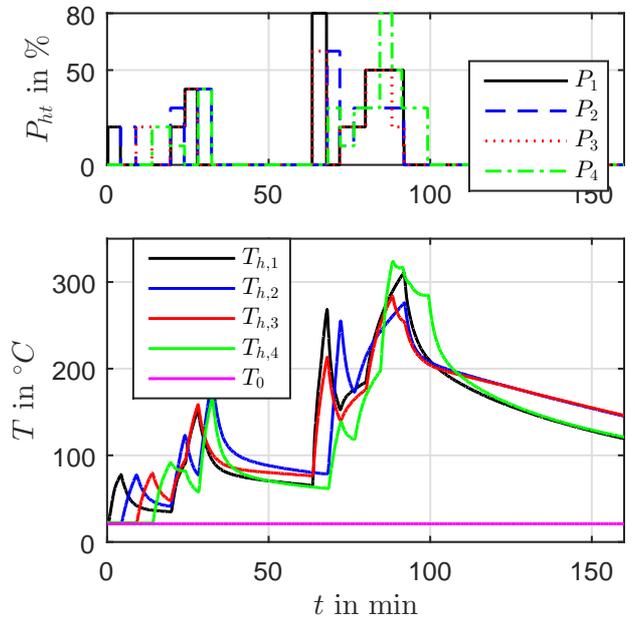}
	\caption{Connection between percentage power $ P $ and the heating tape temperature $ T_{h} $.}
	\label{fig:Vers_THB}
\end{figure}

To verify the accuracy of the created models, the results are compared with measurements. For this purpose, the heating tapes of an industrial extruder in form of\\
Fig. \ref{fig:extruder} (with $ N_{th} = 12 $ and $ N_{h} = 4 $), are driven according to four different electrical power trajectories (in percent) as depicted in the upper half of Fig. \ref{fig:Vers_THB}. The heating tape temperatures are logged by the sensors placed beneath the middle heating tape of each heating zone, representing the input $ \mathbf{u} $ for the mathematical models.\\
We minimize the error between the simulation result and the measurements from the heat up experiment by adapting the unknown heat transfer coefficients between the extruder and the heating tapes $ \alpha_{ht,i} $, between the extruder, and the surrounding $ \alpha_{i} $ and additionally adapting $ c_{p,s} $ and $ \lambda_{s} $. Especially $ \alpha_{ht,i} $ is of great interest, because the heating tape gives off a large amount of the heating power to the environment. A non-linear least square fitting method stated as
\begin{align} \label{eq:nonLinLeastSquare}
	\begin{aligned}
		& \mathbf{p} = \begin{bmatrix}
			c_{p,s} &\lambda_{s} &\alpha_{ht,1} & \cdots &\alpha_{ht,12} & \alpha_{1} & \cdots &\alpha_{n+2m} 
		\end{bmatrix}\\
		& \min _{\mathbf{p}}\|\mathbf{y}_{meas} - \mathbf{y}\|_{2}^{2}\\
		s.t.:\,& \dot{\mathbf{x}} = \mathbf{A}(\mathbf{p})\mathbf{x}+\mathbf{B}(\mathbf{p})\mathbf{u}\\
		&\mathbf{y}=\mathbf{C}_1\mathbf{x}_c
	\end{aligned}
\end{align}
is used. The system starts from a equilibrium state $ \mathbf{x}_0 $ according to the currently measured temperatures. The output $ \mathbf{y} $ include the simulated temperature outcomes at the points where the thermocouples are mounted, therefore $ \mathbf{C}_1 $ does not depend on $ \mathbf{p} $.
The fitting routine \textit{lsqnonlin} from MATLAB is used to calculate $ \mathbf{p} $.\\
\begin{figure*}[tbph]
	\centering
	\includegraphics[scale=0.7]{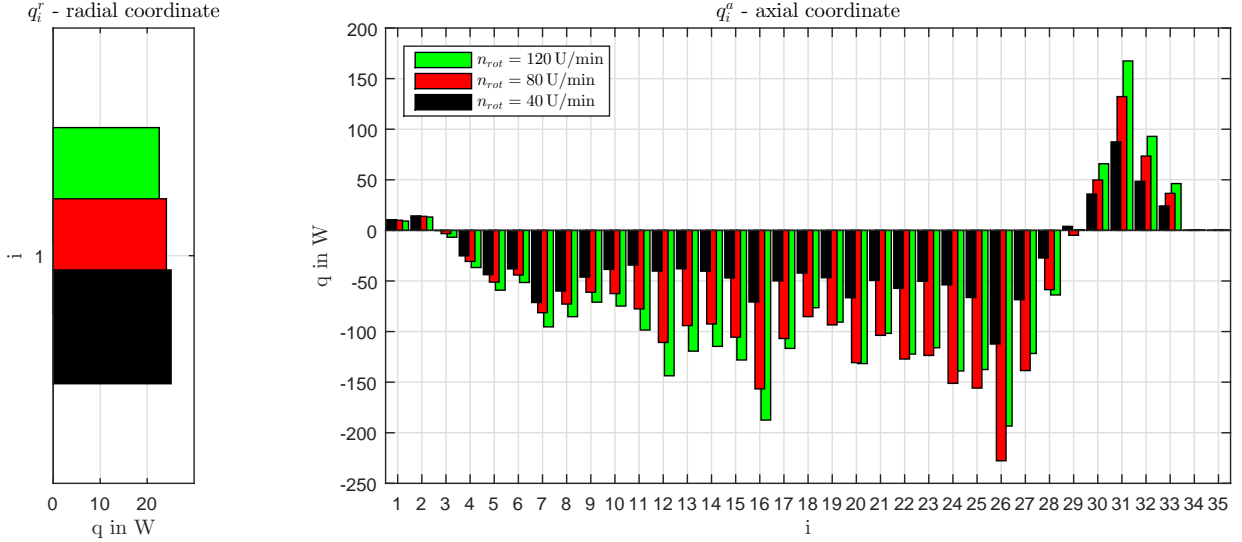}
	\caption{Heat flow observation from experiments by the Smart Sensor, from \emph{Granulate A} in different operating points.}
	\label{fig:stoerN}
\end{figure*}The chosen discretization for the FV model is $ n_x = 35$ elements in axial and $ n_r = 6 $ in radial direction, resulting in a system order of $ N = 210 $. Fig. \ref{fig:Presult} displays the simulation results obtained with the FV (dashed) and the FE (dot-dashed) model, starting from an evenly distributed temperature $ x_{0,i} = T_{0} $, compared to the measurements(solid lines). The small error between the three outcomes confirms a small deviation between the FE and the FV model and also a high accuracy of these models with respect to the measurement. It has to be mentioned that the FV model delivers the result within $ 4 \,\mathrm{s}$ computation time. The FE model comes with a high grid resolution, resulting in a system order of $ 103760 $. A model order reduction was not performed since the result is only intended to serve as a benchmark solution for the FV outcome.
%
%
\section{Smart Sensor}\label{ch:HeatFlowO}
\begin{figure}[H]
	\centering
	\includegraphics[scale=1.0]{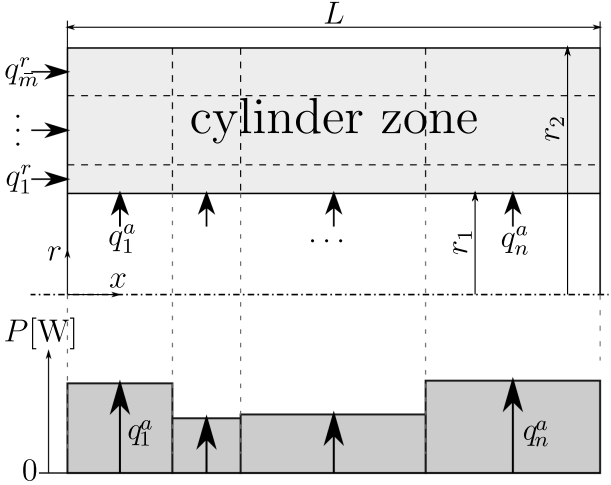}
	\caption{Model adaption as preparation for the Smart Sensor. (e.$ \, $g. $ n=4 $, $ \bar{m}=3 $)}
	\label{fig:cutModel}
\end{figure}

Section \ref{ch:model} proves that the derived FV model reflects the dynamic of the extruder very well. In this Section we design the Smart Sensor based on this model with the aim of observing the thermic behaviour of the material during an extrusion process by estimating the heat flow crossing the surface between the granulate (inside the cylinder) and the extruder cylinder. The FV model is cut along $ \Gamma_{3} $ (Fig. \ref{fig:zones}), resulting in a new model that considers only the cylinder zone, as shown in Fig. \ref{fig:cutModel} depicted as a gray area.
Since the dynamics of the screw conveyor and the granulate are not considered any more we examine the equation \eqref{eq:result} regarding $ \mathbf{x}_c $ 
\begin{align}
	\begin{aligned}
		\dot{\mathbf{x}}_c = \mathbf{A}_{11}\mathbf{x}_c + \mathbf{A}_{12}\mathbf{x}_s + \mathbf{B}_c \mathbf{u}\\
	\end{aligned}\,.
\end{align}
Equation \eqref{eq:connectNode} is used to replace $ \mathbf{A}_{12}\mathbf{x}_s $ by introducing heat flow sources $  \mathbf{q}^a  $ and additionally attach heat sources radially $  \mathbf{q}^r  $  to respect the thermal effect of the actively heated work machine as illustrated in Fig. \ref{fig:cutModel}. The cutting variables, arising by cutting the model as depicted in Fig. \ref{fig:cutModel}, are the heat flows 
\begin{align}\label{eq:cutModelq}
	\begin{aligned}
		\mathbf{q} &= \begin{bmatrix} q_1^a & \cdots & q_n^a, & q_1^r & \cdots & q_{\bar{m}}^r \end{bmatrix}^T\,.
	\end{aligned}
\end{align}
Now $ \mathbf{A}_{12}\mathbf{x}_s $ can be replaced by the cutting variables $ \mathbf{B}_{q}\mathbf{q} $. By respecting $ \mathbf{q} $ as additional state variables the new model results in
\begin{align}\label{eq:cutModelSS1}
	\begin{aligned}
		\dot{\mathbf{x}}_c &= \mathbf{A}_{11}\mathbf{x}_c+\mathbf{B}_q\mathbf{q}+\mathbf{B}_c \mathbf{u}\\
		\dot{\mathbf{q}}&=\mathbf{0}\\
	\mathbf{y}&=\mathbf{C}_1\mathbf{x}_c\,.
		\end{aligned}
\end{align}
The disturbance model reflects the fact that the change rate of the heat flows is sufficiently small, resulting in a high accuracy for the operating states near equilibrium positions. 
The heat flows $ \mathbf{q} $ are interpreted as disturbances acting onto the model and estimate them by the use of the Smart Sensor, which is based on a discretised state space representation of \eqref{eq:cutModelSS1}
\begin{align}\label{eq:cutModelSSE}
	\begin{aligned}
		\begin{bmatrix}
			\mathbf{x}_{c,k+1}\\
			\mathbf{q}_{k+1}
		\end{bmatrix} &= \begin{bmatrix}\bar{\mathbf{A}} & \bar{\mathbf{B}}_{q} \\ \mathbf{0} & \mathbf{E} \end{bmatrix}\begin{bmatrix}
		\mathbf{x}_{c,k}\\
		\mathbf{q}_{k}
	\end{bmatrix}+\bar{\mathbf{B}}_{c}\mathbf{u}_k\\
		\mathbf{y}_k&=\begin{bmatrix} {\mathbf{C}}_1 & \mathbf{0} \end{bmatrix}\begin{bmatrix}
			\mathbf{x}_{c,k}\\
			\mathbf{q}_{k}
		\end{bmatrix}\,.
	\end{aligned}
\end{align}
\begin{figure*}[tbph]
	\centering
	\includegraphics[scale=0.7]{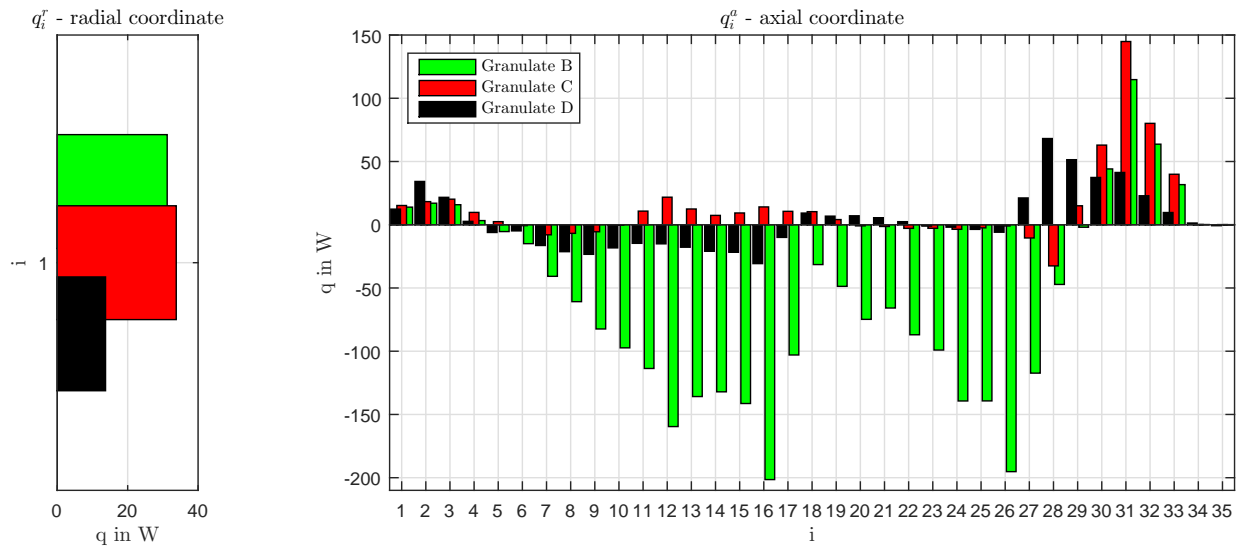}
	\caption{Heat flow observation from extrusion experiments by the Smart Sensor, from different granulate types.}
	\label{fig:stoerVerschK}
\end{figure*}The discretised quantities regarding system \eqref{eq:cutModelSS1} are marked with a bar. Now a LQG-disturbance observer can be designed based on system \eqref{eq:cutModelSSE}. The input for the observer system consists of the heating tape, the ambient and the inner sensors temperatures $ \mathbf{u}_o^T = \begin{bmatrix}\mathbf{u}^T & T_{S_{0}} & \cdots & T_{S_{13}} \end{bmatrix} $. The output is an estimation of the inner temperature distribution $ \mathbf{\hat{x}}_c $ and the heat flows $ \mathbf{\hat{q}}_k $.\\
The stationary observation result coming from an extrusion experiment on a real extruder is depicted in Fig. \ref{fig:stoerN}. The illustrated results are snapshots taken about one hour after extrusion at consistent conditions. Negative heat flows represent a cooling impact due to the extruder cylinder or can be interpreted as an energy consumption effect on the granulate side. The picture on the right side shows the heat flows across the surface between cylinder and granulate, where the filling zone is located around element $ i=35 $ and the extruder end at $ i=1 $. The granulate, mainly existing in a solid state in the area of the filling zone, forces the heating flow $ i=30,\dots,35$ to be positive, because of high friction. The direction of the heat flow changes over the next elements $ i=20,\dots,30$, reflecting a high thermal energy consumption from the granulate. The heat flow decreases to approximately zero at the end of the extruder $ i=0,\dots,20$. The three different colours represent for three different rotational speeds at constant pressure as listed in the legend. The picture on the left shows the observation of the impact of the actively heated work machine mounted to the extruder during the extrusion.\\
The thermal behaviours of different granulate types are illustrates in Fig. \ref{fig:stoerVerschK}. The first material (green) has a strong cooling effect on the cylinder and therefore requires a high heating power to maintain a steady state. The second granulate (black) exhibits a minor cooling effect, the red one even shows a slight heating effect. The observer is capable of observing the thermal impact of the granulates and characterizes them without the knowledge about process parameters of the material.\\
The extruder endowed with a Smart Sensor can be used for granulate behaviour investigations, for observing mechanical wear regarding the screw conveyor or even for optimizing the mechanical extruder build-up itself. The extruder cylinder, combined with the installed temperature sensors, acts like a sensor observing the material behaviour at any time. For this reason, the non-linear material behaviour gets characterized by this sensor without the information of material related parameters and therefore no model of the processed granulate is needed. We are now capable of taking a look inside the extruder, through the help of the developed Smart Sensor. We can use the Smart Sensor as basis for different control concepts to optimize the extrusion process.

\section{CONCLUSIONS}
Two thermal models of a plastic extruder based on the FE and FV approach are derived. We have used the FV model to design a material independent Smart Sensor, which is capable of estimating the energy exchange between the extruder cylinder and the granulate. It changes the descriptive model equations of the respected model, which would be highly non-linear if we add a model of the processed material, to a set of linear system equations. The models are validated by comparing their simulation outcomes to measurements generated by heating experiments on a real extruder. The Smart Sensor has been tested on the real extruder by extrusion experiments. The Smart Sensor can be used as a tool for analysing the granulate behaviour, for checking the work quality from the extruder setup or even for optimizing the mechanical extruder build-up itself. We aim to design a controller, based on the proposed modelling and Smart Sensor strategy, with the goal to control the extruder as fast as possible into a stationary desired extrusion state. 

\begin{ack}
	This work was developed in collaboration with Soplar sa.
	This work has been supported by the FFG, Contract No. 881844: ''Pro$ ^2 $Future''.
\end{ack}
\bibliography{Bibliography}             
                                                   






\end{document}